\documentclass[12pt,a4paper,leqno]{article}

\usepackage{latexsym}
\usepackage{amssymb,exscale}
\usepackage[centertags]{amsmath}
\usepackage{amsthm}
\usepackage[all]{xy}
\usepackage{graphicx}

\usepackage[dvips]{hyperref}

\numberwithin{equation}{section}

\swapnumbers
\theoremstyle{definition}
\newtheorem{theorem}[equation]{Theorem}
\newtheorem*{thm}{Theorem}
\newtheorem{lemma}[equation]{Lemma}
\newtheorem{proposition}[equation]{Proposition}
\newtheorem{corollary}[equation]{Corollary}
\newtheorem{definition}[equation]{Definition}
\newtheorem{remark}[equation]{Remark}
\newtheorem{question}[equation]{Question}

\renewcommand{\phi}{\varphi}

\newcommand{\D}{\mathrm{d}}

\renewcommand{\(}{\bigl(}
\renewcommand{\)}{\bigr)\vphantom{)}}

\newcommand{\ip}[2]{\langle#1,#2\rangle}

\newcommand{\imply}{\;\;\,\Longrightarrow\;\;\,}
\newcommand{\imp}{$ \Longrightarrow $ }

\newcommand{\DCS}{\operatorname{DCS}}
\newcommand{\Exp}{\operatorname{Exp}}
\newcommand{\mes}{\operatorname{mes}}

\newcommand{\modO}{{\operatorname{mod}\,0}}

\newcommand{\One}{\mathbf1}

\newcommand{\eps}{\varepsilon}
\newcommand{\si}{\sigma}

\newcommand{\om}{\omega}
\newcommand{\Om}{\Omega}
\newcommand{\de}{\delta}
\newcommand{\De}{\Delta}

\newcommand{\Ec}{\mathcal E}
\newcommand{\F}{\mathcal F}
\newcommand{\G}{\mathcal G}
\newcommand{\A}{\mathcal A}

\newcommand{\const}{{\mathrm{const}}}
\newcommand{\Ex}{\mathbb E\,}
\renewcommand{\Pr}[1]{\mathbb{P}\mskip1.5mu\(\mskip1.5mu#1\mskip1.5mu\)}
\newcommand{\R}{\mathbb R}
\newcommand{\Q}{\mathbb Q}

\newcommand{\cE}[2]{\mathbb{E}\mskip1.5mu\(\mskip1.5mu#1\mskip1.5mu
 \big|\mskip1.5mu#2\mskip1.5mu\)}

\newcommand{\cP}[2]{\mathbb{P}\mskip1.5mu\(\mskip1.5mu#1\mskip1.5mu
 \big|\mskip1.5mu#2\mskip1.5mu\)}

\newcommand{\sif}{$\sigma$\nobreakdash-field}

\newcommand{\my}[2]{#1\!\times\!B\,\cup\,B\!\times\!#2}

\newcommand{\almost}[1]{$#1$\nobreakdash-\hspace{0pt}almost}

\newcommand{\valued}[1]{$#1$\nobreakdash-\hspace{0pt}valued}

\newcommand{\saturated}[1]{$#1$\nobreakdash-\hspace{0pt}saturated}

\begin{document}

\title{Brownian local minima and other\\ random dense countable sets} 

\author{Boris Tsirelson}

\date{}
\maketitle

\begin{abstract}
We compare two examples of random dense countable sets, \emph{Brownian local
minima} and \emph{unordered uniform infinite sample}. They appear to be
identically distributed. A framework for such notions is proposed. In
addition, random elements of other singular spaces (especially, reals modulo
rationals) are considered.
\end{abstract}

\section*{Introduction}
For almost every Brownian path $ \om = (b_t)_{t\in[0,1]} $ on $ [0,1] $, the
set
\begin{equation}\label{0.1}
M_\om = \{ s \in (0,1) : \exists \eps>0 \;\> \forall t \in (s-\eps,s) \cup
(s,s+\eps) \;\> b_s < b_t \}
\end{equation}
of local minimizers on $ (0,1) $ is a dense countable subset of $ (0,1)
$. Should we say that $ (M_\om)_\om $ is a random countable dense set? Can we
give an example of an event of the form $ \{ \om : M_\om \in A \} $ possessing
a probability different from $ 0 $ and $ 1 \, $? No, we cannot (see also
Corollary \ref{5.1}). All dense
countable subsets of $ (0,1) $ are a set $ \DCS(0,1) $ (of sets), just a set,
not a  Polish space, not even a standard Borel space. What should we mean by
an \valued{\DCS(0,1)} random variable and its distribution? Apart from such
conceptual questions we have specific examples and questions; here is one. A
`uniform infinite sample', that is, an infinite sample from the uniform
distribution on $ (0,1) $ may be described by the product $ \( (0,1)^\infty,
\mes^\infty \) $ of an infinite sequence of copies of the probability space $
\( (0,1), \mes \) $, where `$ \mes $' stands for the Lebesgue measure on $
(0,1) $. For almost every point $ u = (u_1,u_2,\dots) $ of this product space,
the set
\begin{equation}\label{0.2}
S_u = \{ u_1, u_2, \dots \} = \{ s \in (0,1) : \exists n \;\> u_n = s \}
\end{equation}
is a dense countable subset of $ (0,1) $. It appears that $ (M_\om)_\om $ and
$ (S_u)_u $ are identically distributed in the following sense (see
Th.~\ref{4.7}).

\begin{thm}
There exists a joining $ J $ between the Brownian motion on $ [0,1] $ and the
uniform infinite sample such that $ M_\om = S_u $ for \almost{J} all pairs $
(\om,u) $.
\end{thm}

The theorem follows from a more general theory presented below. If you
consider the theory too general, try to find a better proof of this theorem or
maybe its two-point corollary; namely, construct (at least) two independent
uniform random variables $ U_1, U_2 $ coupled with the Brownian motion $
(B_t)_t $ in such a way that almost surely $ U_1, U_2 $ are (some of the)
local minimizers of $ (B_t)_t $.

\section[]{\raggedright Definitions}
\label{sect1}The set $ \DCS(0,1) $ is a \emph{singular space} in the sense of Kechris
\cite[\S2]{Ke99}: a `bad' quotient space of a `good' space by a `good'
equivalence relation. (A simpler example of a singular space is $ \R / \Q $,
reals modulo rationals.) Namely,
\begin{equation}\label{1.1}
\DCS(0,1) = (0,1)^\infty_{\ne} / E \, .
\end{equation}
Here $ (0,1)^\infty_{\ne} $ is the set of all sequences $ u = (u_1,u_2,\dots)
$ of pairwise different points of $ (0,1) $, and $ E $ is the following
equivalence relation on $ (0,1)^\infty_{\ne} $:
\begin{equation}\label{1.2}
E = \{ (u,v) : S_u = S_v \} \, ,
\end{equation}
$ S_u $ being defined by \eqref{0.2}. (In fact, equivalence classes are orbits
of a natural action of the infinite
permutation group, see \cite[Sect.~2e]{Ts04}.) Note that $ (0,1)^\infty_{\ne}
$ is a standard Borel space and $ E $ is a Borel subset of $ (0,1)^\infty_{\ne}
\times (0,1)^\infty_{\ne} $. It is possible to equip the quotient space with
its natural \sif\ (of sets whose inverse images are measurable) and define
random variables and distributions accordingly. Is it a good idea? I do not
know. (See also Sect.~\ref{sect5}.) I prefer another concept of a random
element in a singular space, sketched in \cite[Sect.~2e]{Ts04} and formalized
below.

\emph{Throughout Sections \textup{1--4}, either by assumption or by
construction, all probability spaces are standard.} Recall that a standard
probability space (known also as a Lebesgue-Rokhlin space) is a probability
space isomorphic $ (\modO) $ to an interval with the Lebesgue measure, a
finite or countable collection of atoms, or a combination of both.

\begin{definition}\label{1.3}
Let $ B $ be a standard Borel space, $ E \subset B \times B $ an equivalence
relation on $ B $, and $ \Om $ (or rather $ (\Om,\F,P) $) a probability
space. A map $ X : \Om \to B/E $ is called \emph{measurable,} if there exists
a measurable map $ Y : \Om \to B $ such that the following diagram is
commutative:
\[
\xymatrix{
 \Om \ar[r]^{Y} \ar[dr]_{X} & B \ar[d]^{\text{canonical projection}}
\\
 & B/E
}
\]
\end{definition}

Note that $ B $ and $ E $ have to be given. We do not touch on the question,
what happens if (in some sense) $ B/E = B_1/E_1 $. Note also that
Def.~\ref{1.3} is in the spirit of the `diffeology' (see \cite{Ig}, especially
Sect.~1.14 `Quotient of manifolds' and 1.15 `The irrational torus').

Equivalence classes of measurable maps $ \Om \to B $ are elements of the set $
L_0 ( \Om \to B ) $ of \valued{B} random variables on $ \Om $. Similarly, we
define $ L_0 ( \Om \to B/E ) $ as the set of all equivalence classes of
measurable maps $ \Om \to B/E $ (the equivalence being the equality almost
everywhere, as usual). Being equipped with the natural \sif, the set $ L_0 (
\Om \to B ) $ is a standard Borel space. The set $ L_0 ( \Om \to B/E ) $ may
be treated as a singular space,
\[
L_0 ( \Om \to B/E ) = L_0 ( \Om \to B ) / L_0 ( \Om \to E ) \, ,
\]
where $ L_0 ( \Om \to E ) $ is the following equivalence relation on $ L_0 (
\Om \to B ) $: $ (f,g) \in L_0 ( \Om \to E ) $ iff $ (f(\om),g(\om)) \in E $
for almost all $ \om $.

Def.~\ref{1.3} is compatible with the usual definition in the following
sense. Let $ A,B $ be two standard Borel spaces, $ f : B \to A $ a Borel
function, and $ E = \{ (x_1,x_2) : f(x_1)=f(x_2) \} $. Then $ B/E = A $ (after
the evident identification). It is easy to check that
\[
L_0 ( \Om \to B/E ) = L_0 ( \Om \to A )
\]
(after the evident identification); here $ L_0 ( \Om \to A ) $ is defined as
usual, while $ L_0 ( \Om \to B/E ) $ is defined by \ref{1.3}.

Waiving the \sif\ on $ B/E $ we lose the usual definition of a distribution on
$ B/E $. Instead we may define the notion `identically distributed' as
follows.

\begin{definition}\label{1.4}
Let $ B $ be a standard Borel space, $ E \subset B \times B $ an equivalence
relation on $ B $, and $ \Om_1, \Om_2 $ probability spaces. Random variables $
f \in L_0 ( \Om_1 \to B/E ) $, $ g \in L_0 ( \Om_2 \to B/E ) $ are
\emph{identically distributed,} if there exist a probability space $ \Om $ and
measure preserving maps $ T_1 : \Om \to \Om_1 $, $ T_2 : \Om \to \Om_2 $ such
that $ f(T_1(\om)) = g(T_2(\om)) $ for almost all $ \om \in \Om $.
\end{definition}

That is, the following diagram must be commutative ($ \modO $):
\[
\xymatrix{
 & \Om \ar[dl]_{T_1} \ar[dr]^{T_2}
\\
 \Om_1 \ar[dr]_f & & \Om_2 \ar[dl]^g
\\
 & B/E
}
\]
The joint distribution of $ T_1(\om), T_2(\om) $ is a joining, that is, a
measure $ J $ on $ \Om_1 \times \Om_2 $ with given marginals $ P_1, P_2
$. Here is a definition equivalent to \ref{1.4}: $ f,g $ are identically
distributed, if
there exists a joining $ J $ between $ \Om_1 $ and $ \Om_2 $ such that $
f(\om_1) = g(\om_2) $ for \almost{J} all pairs $ (\om_1,\om_2) $.

Def.~\ref{1.4} is compatible with the usual definition, similarly to
Def.~\ref{1.3}. Namely, let $ B/E = A $ be a standard Borel space. Then $ f
\in L_0 ( \Om_1 \to B/E ) $, $ g \in L_0 ( \Om_2 \to B/E ) $ are identically
distributed according to Def.~\ref{1.4} if and only if $ f \in L_0 ( \Om_1 \to
A ) $, $ g \in L_0 ( \Om_2 \to A ) $ are identically distributed in the usual
sense.

\begin{definition}\label{1.5}
Let $ B $ be a standard Borel space and $ E \subset B \times B $ an
equivalence relation on $ B $. A \emph{distribution} on $ B/E $ is an
equivalence class of \valued{B/E} random variables on $ \( (0,1), \mes \) $;
here equivalence of two random variables means that they are identically
distributed.
\end{definition}

Def.~\ref{1.5} is compatible with the usual definition (similarly to
\ref{1.3}, \ref{1.4}).

We return to $ \DCS(0,1) $ treated as $ B/E $ according to \eqref{1.1},
\eqref{1.2}. The first example of a \valued{\DCS(0,1)} random variable is the
\emph{unordered uniform infinite sample.} We define it as the \valued{B/E}
random variable corresponding to the ordered uniform infinite sample. The
latter is the \valued{B} random variable $ (U_1,U_2,\dots) $; here $ B =
(0,1)^\infty_{\ne} $ and $ U_1,U_2,\dots $ are i.i.d.\ random variables
uniform on $ (0,1) $. The unordered uniform infinite sample depends on the
choice of $ U_1,U_2,\dots $ and the underlying probability space, but its
distribution is uniquely determined.

\section[]{\raggedright Main lemma}
\label{sect2}\begin{lemma}\label{2.1}
Let $ X_1,X_2,\dots $ be real-valued random variables (on some probability
space) such that for every $ n = 1,2,\dots $ the conditional distribution of $
X_n $ given $ X_1,\dots,X_{n-1} $ has a density $ (x,\om) \mapsto f_n(x,\om)
$. If
\[
\sum_{n=1}^\infty f_n(x,\om) = \begin{cases}
 \infty & \text{for } 0<x<1,\\
 0 & \text{otherwise}
\end{cases}
\]
for almost all $ x $ and $ \om $, then the \valued{\DCS(0,1)} random variable
\[
\om \mapsto \{ X_1(\om), X_2(\om), \dots \} = \{ x\in\R : \exists n \;\>
X_n(\om)=x \}
\]
is distributed like an unordered uniform infinite sample.
\end{lemma}

The proof is given below after some discussion.
We see that the distribution of an unordered (not just uniform) infinite
sample does not depend on the underlying one-dimensional distribution on $
(0,1) $ provided that the latter distribution has a strictly positive density
on $ (0,1) $. The same holds for independent (not just identically
distributed) $ X_n $, provided that each $ X_n $ has a density $ f_n $ and $
f_1 + f_2 + \dots = \infty $ almost everywhere on $ (0,1) $. Especially,
the case $ f_1 = f_3 = \dots $, $ f_2 = f_4 = \dots $ leads to the following
fact.

\begin{corollary}
If $ \Om_1 \ni \om_1 \mapsto A(\om_1) \in \DCS(a,b) $ is an unordered uniform
infinite sample on $ (a,b) $ and $ \Om_2 \ni \om_2 \mapsto B(\om_2) \in
\DCS(b,c) $ is an unordered uniform infinite sample on $ (b,c) $ then 
\[
\Om_1 \times \Om_2 \ni (\om_1,\om_2) \mapsto A(\om_1) \cup B(\om_2) \in
\DCS(a,c) 
\]
is distributed like an unordered uniform infinite sample on $ (a,c) $. 
\end{corollary}

\begin{proof}[Proof of Lemma \textup{\ref{2.1}}]
We introduce a Poisson random subset of the strip $ (0,1) \times (0,\infty) $
on some probability space $ \Om $,
\[
\Om \ni \om \mapsto A(\om) \subset (0,1) \times (0,\infty) \, ,
\]
whose intensity measure is the (two-dimensional) Lebesgue measure on the
strip. Almost surely, $ A(\om) $ is a countable, locally finite set. We define
functions $ g_n : (0,1)^n \to [0,\infty) $ by
\[
f_n (x,\om) = g_n \( X_1(\om), \dots, X_{n-1}(\om), x \)
\]
(some ambiguity in $ g_n $ is harmless) and construct random variables $
Y_1, Y_2, \dots : \Om \to (0,1) $ and $ T_1, T_2, \dots : \Om \to (0,\infty)
$ step by step, as follows.

The first step:
\[
T_1(\om) = \min \{ t>0 : \exists y \in (0,1) \;\> (y,tg_1(y)) \in A(\om) \} =
\min_{(y,h)\in A(\om)} \frac{ h }{ g_1(y) } \, ;
\]
this random variable is distributed $ \Exp(1) $, since $ \int_0^1 g_1(y) \, \D
y = 1 $. The corresponding point $ y $ (evidently unique a.s.) gives us $
Y_1(\om) $,
\[
\( Y_1(\om), T_1(\om) g_1(Y_1(\om)) \) \in A(\om) \, .
\]
\[
\begin{gathered}\includegraphics{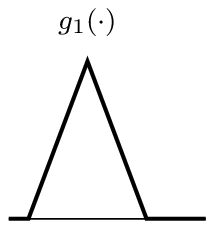}
 \qquad\qquad\includegraphics{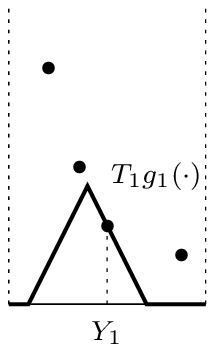}\end{gathered}
\]
The random variable $ Y_1 $ is distributed like $ X_1 $ (since $ g_1 $ is its
density) and independent of $ T_1 $.

Probabilistic statements about the second step (below) are conditioned on $
T_1 $ and $ Y_1 $. The conditioning does not perturb the Poisson set $ A $
above the graph of the function $ T_1 g_1(\cdot) $.

The second step:
\[
T_2(\om) = \min \{ t>0 : \exists y \in (0,1) \;\> (y,T_1(\om)g_1(y) +
t g_2(Y_1(\om),y)) \in A(\om) \}
\]
is distributed $ \Exp(1) $ (since $ \int_0^1 g_2 (Y_1(\om),y)) \, \D y = 1 $
a.s.), and we define $ Y_2 $ as the unique $ y $,
\[
\( Y_2(\om), T_1(\om) g_1(Y_2(\om)) + T_2(\om) g_2(Y_1(\om),Y_2(\om)) \) \in
A(\om) \, .
\]
\[
\begin{gathered}\includegraphics{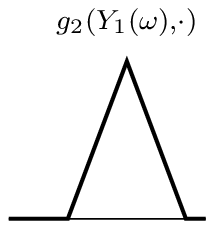}
 \qquad\qquad\includegraphics{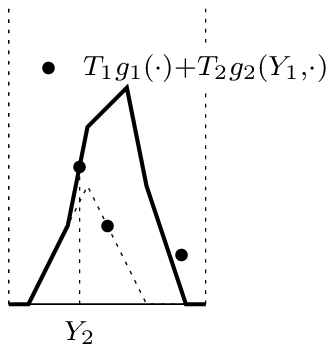}\end{gathered}
\]
Random variables $ T_2, Y_2 $ are independent; $ T_2 $ is distributed $
\Exp(1) $, while $ Y_2 $ has the density $ g_2(Y_1,\cdot) $. These relations
are conditional; unconditionally, the pair $ (Y_1,Y_2) $ is distributed like $
(X_1,X_2) $ and independent of the pair $ (T_1,T_2) \sim \Exp(1) \otimes
\Exp(1) $. Conditioning on $ T_1, Y_1, T_2, Y_2 $ does not perturb the
Poisson set $ A $ above the graph of the function $ T_1 g_1(\cdot) + T_2
g_2(Y_1,\cdot) $.

Continuing the process we get random variables $ T_n, Y_n $ ($ n=1,2,\dots $)
on $ \Om $ such that the sequence $ (Y_1,Y_2,\dots) $ is distributed like $
(X_1,X_2,\dots) $ and independent of the i.i.d.\ sequence $ (T_1,T_2,\dots) $
of $ \Exp(1) $ random variables. Conditioning on all $ Y_n $ and $ T_n $ does
not perturb the Poisson set $ A $ above the graph of the function $ \sum_n
T_n g_n(Y_1,\dots,Y_{n-1},\cdot) $ (a void claim if the sum is infinite
everywhere).

Now we use the condition $ \sum_n f_n(x,\om) = \infty $. It gives us
\[
\sum_n g_n ( Y_1(\om),\dots,Y_{n-1}(\om), x ) = \infty
\]
for almost all $ \om \in \Om $, $ x \in (0,1) $. It follows that
\[
\sum_n T_n(\om) g_n ( Y_1(\om),\dots,Y_{n-1}(\om), x ) = \infty
\]
for almost all $ \om \in \Om $, $ x \in (0,1) $ (since the relation holds
conditionally, given $ (Y_1,Y_2,\dots) $). We see that almost no points of the
strip remain above the graph of this sum, and therefore, no one point of $
A(\om) $ does (a.s.). All points of $ A(\om) $ are used in our
construction. Therefore the \valued{\DCS(0,1)} random variable $ \om \mapsto
\{ Y_(\om), Y_2(\om), \dots \} $ is just the projection of the Poisson set $
A(\om) $, therefore, an unordered uniform infinite sample. On the other hand,
$ \{ Y, Y_2, \dots \} $ is distributed like $ \{ X, X_2, \dots \} $.
\end{proof}

\section[]{\raggedright A sufficient condition}
\label{sect3}The condition $ \sum f_n(x,\om) = \infty $ of Lemma \ref{2.1} may be checked
pointwise. For every $ x $ we have a series of random variables $ f_n(x,\cdot)
\ge 0 $, and check its divergence a.s. If this holds for all (or almost all) $
x \in (0,1) $, Lemma \ref{2.1} is applicable.

Let $ Y_1, Y_2, \dots : \Om \to [0,\infty) $ be a sequence of random variables
(generally, interdependent). We seek a sufficient condition for the property
\begin{equation}\label{3.1}
\sum Y_n = \infty \quad \text{a.s.}
\end{equation}
If $ \sum Y_n < \infty $ a.s.\ then $ Y_n \to 0 $ a.s., which implies $ \Pr{
Y_n < \eps } \to 1 $ for any $ \eps > 0 $ (since indicators $
\One_{[\eps,\infty)} (Y_n) $ converge to $ 0 $ a.s.). Given an event $ A
\subset \Om $, $ \Pr{A} > 0 $, we may apply the remark above to the
probability space $ A $ (with the conditional measure). The case $ A = \{ \sum
Y_n < \infty \} $ leads to $ \liminf_n \Pr{Y_n<\eps} \ge \Pr{ \sum Y_n <
\infty } $, thus,
\[
\Pr{ \sum Y_n < \infty } \le \lim_{\eps\to0+} \liminf_n \Pr{Y_n<\eps} \le
\lim_{\eps\to0+} \sup_n \Pr{Y_n<\eps} \, .
\]
We see that the condition
\[
\lim_{\eps\to0+} \sup_n \Pr{Y_n<\eps} = 0
\]
is sufficient for \eqref{3.1}. Unfortunately, this sufficient condition is too
strong for our purpose. We assume a weaker condition
\begin{equation}\label{3.2}
\lim_{\eps\to0+} \sup_n \Pr{0<Y_n<\eps} = 0 \, .
\end{equation}
Surely, \eqref{3.2} does not imply \eqref{3.1}, since $ Y_n $ may vanish. We
seek an additional condition on the events $ \{ Y_n = 0 \} $.

Once again, if $ \sum Y_n < \infty $ a.s.\ then $ \Pr{ Y_n < \eps } \to 1 $
for any $ \eps $; combined with \eqref{3.2} it gives $ \Pr{Y_n=0} \to 1 $. As
before, we condition on the event $ \{ \sum Y_n < \infty \} $ (which does not
invalidate \eqref{3.2}; of course we assume here that the event is of positive
probability). The straightforward conclusion
\[
\liminf_n \Pr{Y_n=0} \ge \Pr{ \sum Y_n < \infty }
\]
is of little interest; instead, we introduce the condition
\begin{equation}\label{3.3}
\Pr{ A \cap \{Y_n>0\} } \to 0 \quad \text{implies} \quad \Pr A = 0
\end{equation}
for all measurable sets $ A \subset \Om $.

\begin{lemma}\label{3.4}
Every sequence $ (Y_n) $ satisfying both \eqref{3.2} and \eqref{3.3} satisfies
\eqref{3.1}.
\end{lemma}

\begin{proof}
Otherwise $ \cP{ Y_n=0 }{ A } \to 1 $ where $ A = \{ \sum Y_n < \infty \} $, $
\Pr A > 0 $. Thus, $ \cP{ Y_n>0 }{ A } \to 0 $ in contradiction to
\eqref{3.3}.
\end{proof}

Now we need a condition sufficient for \eqref{3.3}. Let $ T : \Om \to \Om $ be
a (strongly) mixing measure preserving transformation and $ B \subset \Om $ a
measurable set, $ \Pr B > 0 $. Then
\[
\Pr{ A \cap T^{-n}(B) } \to \Pr A \Pr B
\]
for every measurable set $ A \subset \Om $, which ensures \eqref{3.3} if the
events $ \{ Y_n > 0 \} $ are of the form $ T^{-n}(B) $. However, we need a
more general case,
\[
\{ Y_n > 0 \} = T^{-n}(A_n)
\]
where $ \{ A_1, A_2, \dots \} $ is a precompact set (of events). The
precompactness means that every subsequence $ (A_{n_k})_k $ contains a
subsequence $ (A_{n_{k_i}})_i $ such that $ \Pr{ A_{n_{k_i}} \setminus
A_{n_{k_j}} } \to 0 $ as $ i,j\to\infty $. Or equivalently, all indicator
functions $ \One_{A_n} $ belong to a single compact subset of $ L_2(\Om,\F,P)
$.

\begin{lemma}\label{3.5}
Let $ T : \Om \to \Om $ be a mixing measure preserving transformation and $
A_n \subset \Om $ measurable sets such that
\[
\limsup_n \Pr{ A_n } > 0
\]
and $ \{ A_1, A_2, \dots \} $ is a precompact set. Then $ \Pr{ A \cap
T^{-n}(A_n) } \to 0 $ implies $ \Pr A = 0 $ for all measurable sets $ A
\subset \Om $.
\end{lemma}

\begin{proof}
The isometric operator $ U : L_2 \to L_2 $ defined by $ (Uf)(\om) = f(T\om) $
satisfies
\[
U^n \to \Ex \quad \text{weakly as $ n\to\infty $}
\]
where $ \Ex $ is the expectation treated as the projection onto the
one-dimensional space of constants. It follows that for every $ g \in L_2 $
the convergence
\[
\ip{ U^n f }{ g } - \ip{ \Ex f }{ g } \to 0
\]
is uniform in $ f $ as long as $ f $ runs over a compact set. We take $ f =
f_n = \One_{A_n} $, $ g = \One_A $ and get
\[
\Pr{ T^{-n}(A_n) \cap A } - \Pr{A_n} \Pr{A} = \ip{ U^n f_n }{ g } - \ip{ \Ex
f_n }{ g } \to 0 \, .
\]
If $ A $ satisfies $ \Pr{ A \cap T^{-n}(A_n) } \to 0 $ then $ \Pr{A_n} \Pr{A}
\to 0 $ which implies $ \Pr A = 0 $.
\end{proof}

The following proposition combines the ideas of \ref{3.4}, \ref{3.5} and
introduces one more idea (the transition from $ Z_n $ to $ Y_n $) needed for
the next section.

\begin{proposition}\label{3.6}
Let $ (\Om,\F,P) $ be a probability space, $ T : \Om \to \Om $ a mixing
measure preserving transformation, $ Y_1, Y_2, \dots : \Om \to [0,\infty) $
and $ Z_1, Z_2, \dots : \Om \to [0,\infty) $ random variables, and $ A_1, A_2,
\dots \subset \Om $ a precompact sequence of measurable sets such that $
\limsup_n \Pr{A_n} > 0 $. Assume that

(a) $ \{ Y_n \ne 0 \} = \{ Z_n \ne 0 \} = T^{-n} (A_n) $ for each $ n $,

(b) $ Y_n = \cE{ Z_n }{ Y_n } $ for each $ n $,

(c) $ \lim_{\eps\to0+} \sup_n \Pr{0<Z_n<\eps} = 0 $.

Then $ \sum Y_n = \infty $ a.s.
\end{proposition}

\begin{proof}
First, we claim that
\begin{equation}\label{3.7}
\Pr{ 0 < Y_n < \eps } \le 2 \Pr{ 0 < Z_n < 2\eps }
\end{equation}
for all $ n $ and $ \eps $. Proof: conditioning on the event $ \{ Y_n \ne 0 \}
= \{ Z_n \ne 0 \} $ reduces \eqref{3.7} to a simpler claim: $ \Pr{ Y_n < \eps
} \le 2 \Pr{ Z_n < 2\eps } $ for any random variables $ Y,Z : \Om \to
[0,\infty) $ such that $ Y = \cE Z Y $. We note that
\[
\cP{ Z \ge 2\eps }{ Y } \le \frac1{2\eps} \cE{ Z }{ Y } = \frac1{2\eps} Y \, ,
\]
thus,
\begin{multline*}
\Pr{ Y < \eps } = \Pr{ \tfrac1{2\eps} Y < \tfrac12 } \le \Pr{ \cP{ Z \ge
 2\eps }{ Y } < \tfrac12 } = \\
= \Pr{ \cP{ Z < 2\eps }{ Y } > \tfrac12 } \le \( \tfrac12 \)^{-1}
 \Ex \Big( \cP{ Z < 2\eps }{ Y } \Big) = 2 \Pr{ Z < 2\eps } \, ,
\end{multline*}
which proves the claim.

Combining \eqref{3.7} with (c) we see that the sequence $ (Y_n)_n $ satisfies
\eqref{3.2}. On the other hand, Lemma \ref{3.5} combined with (a) gives
\eqref{3.3}. Lemma \ref{3.4} completes the proof.
\end{proof}

\begin{remark}\label{3.8}
Condition \ref{3.6}(a) may be relaxed:
$ \{ Y_k \ne 0 \} = \{ Z_k \ne 0 \} = T^{-n_k} (A_k) $ for some $ n_1 < n_2 <
\dots $
\end{remark}

\section[]{\raggedright Main theorem}
\label{sect4}We consider the usual one-dimensional Brownian motion $ (t,\om) \mapsto
B_t(\om) $ for $ t \in [0,1] $; $ \om $ runs over a probability space $
(\Om,\F,P) $. The set $ M_\om $ of all local minimizers of the path $ t
\mapsto B_t(\om) $ on $ (0,1) $ is well-known to be a dense countable set,
\[
M_\om \in \DCS(0,1) \quad \text{for almost all $ \om $.}
\]

\begin{lemma}\label{4.1}
The map $ \om \to M_\om $ from $ \Om $ to $ \DCS(0,1) $ is measurable (as
defined by \ref{1.3} using \eqref{1.1}).
\end{lemma}

\begin{proof}
We need a measurable enumeration of $ M_\om $, that is, a sequence of random
variables $ X_1, X_2, \dots : \Om \to (0,1) $ such that
\begin{equation}\label{4.2}
\begin{gathered}
M_\om = \{ X_1(\om), X_2(\om), \dots \} \, , \\
X_k(\om) \ne X_l(\om) \quad \text{for $ k \ne l $}
\end{gathered}
\end{equation}
for almost all $ \om $. We enumerate all dyadic intervals by the numbers $
2,3,4,\dots $,
\[
\begin{matrix}
n\;\; & 2 & 3 & 4 & 5 & 6 & 7 & 8 & 9 & \dots \\
I_n\;\; & (0,1) & (0,\frac12) & (\frac12,1) & (0,\frac14) & (\frac14,\frac12) &
(\frac12,\frac34) & (\frac34,1) & (0,\frac18) & \dots
\end{matrix}
\]
For each $ n > 1 $ we consider the left half $ I'_n = I_{2n-1} $ and the right
half $ I''_n = I_{2n} $ of $ I_n $, the corresponding Brownian minimizers $
X'_n, X''_n $,
\[
X'_n \in I'_n, \;\; B_{X'_n} = \inf_{t\in I'_n} B_t, \quad
X''_n \in I''_n, \;\; B_{X''_n} = \inf_{t\in I''_n} B_t,
\]
and define $ X_n $ as the minimizer that corresponds to the greater minimum,
\[
X_n = \begin{cases}
 X'_n &\text{if $ B_{X'_n} > B_{X''_n} $},\\
 X''_n &\text{if $ B_{X'_n} < B_{X''_n} $}.
\end{cases}
\]
In addition we define $ X_1 $ as the Brownian minimizer on the whole $ (0,1)
$.

For every $ k=0,1,2,\dots $ the $ 2^k $ numbers $ X_1,\dots,X_{2^k} $ are
nothing but the Brownian minimizers on the $ 2^k $ dyadic intervals $ I_n $
for $ 2^k < n \le 2^{k+1} $, that is, the intervals $ \( (i-1)/2^k, i/2^k \) $
for $ i=1,\dots,2^k $ (randomly rearranged, of course). Therefore \eqref{4.2}
is satisfied.
\end{proof}

\begin{lemma}\label{4.3}
The random variables $ X_1, X_2, \dots $ introduced in the proof of Lemma
\ref{4.1} are such that for every $ n=1,2,\dots $ the conditional distribution
of $ X_n $ given $ X_1,\dots,X_{n-1} $ has a density $ (x,\om) \mapsto
f_n(x,\om) $.
\end{lemma}

\begin{proof}
We define by $ \Ec_n $ the sub-\sif\ of $ \F $ generated by $
X_1,\dots,X_{n-1} $, and by $ C_n $ the event $ \{ X_n \in I'_n \} $. Note
that $ C_n \in \Ec_{n-1} $, since $ C_n = \{ X_1 \in I''_n \} \cup \dots \cup
\{ X_{n-1} \in I''_n \} $. Note also that $ C_n = \{ X_n \in I_{2n-1} \} $ and
$ \Om \setminus C_n = \{ X_n \in I_{2n} \} $. We define by $ \G_n $ the
sub-\sif\ of $ \F $ generated by all $ B_s $ for $ s \in [0,1] \setminus I_n
$, and by $ \F_n $ the sub-\sif\ of $ \F $ that contains $ C_n $, coincides
with $ \G_{2n-1} $ on $ C_n $ and with $ \G_{2n} $ on $ \Om \setminus C_n
$. In other words, $ \F_n $ consists of sets of the form $ ( A \cap \{ X_n \in
I_{2n-1} \} ) \cup ( B \cap \{ X_n \in I_{2n} \} ) $ for $ A \in \G_{2n-1} $,
$ B \in \G_{2n} $.

We claim that $ \Ec_{n-1} \subset \F_n $. Proof: both \sif s contain $ C_n $;
on $ C_n $ the inclusion holds since here $ X_1,\dots,X_{n-1} \in [0,1]
\setminus I'_n $; on $ \Om \setminus C_n $ the inclusion holds since here $
X_1,\dots,X_{n-1} \in [0,1] \setminus I''_n $.

The conditional distribution of $ X_n $ given $ \F_n $ is easy to describe. On
$ C_n $ it is the conditional distribution of the Brownian minimizer on $ I'_n
$ under three conditions. Two conditions are boundary values of the Brownian
path on the two endpoints of $ I'_n $. The third condition is a lower
bound on (the minimum of) the Brownian path on $ I'_n $; it must exceed the
minimum on $ I''_n $. A similar description holds on $ \Om \setminus C_n
$. Clearly, the conditional distribution of $ X_n $ given $ \F_n $ has a
density $ (x,\om) \mapsto g_n(x,\om) $ (see also \eqref{4.4} below).

Taking into account that $ \Ec_{n-1} \subset \F_n $ we conclude that the
conditional distribution of $ X_n $ given $ \Ec_{n-1} $ has a density $
(x,\om) \mapsto f_n(x,\om) $,
\[
f_n (x,\cdot) = \cE{ g_n(x,\cdot) }{ \Ec_{n-1} } \, .
\]
\end{proof}

Here is an explicit formula for the conditional density $ g_n $ introduced
above: for $ \om \in C_n $ and $ x \in I'_n = (u,v) $,
\begin{multline}\label{4.4}
g_n(x,\om) = \\
\frac1{v-u} \phi \bigg( \frac1{\sqrt{v-u}} \Big( B_u(\om) -
 \min_{I''_n} B(\cdot,\om) \Big), \frac1{\sqrt{v-u}} \Big( B_v(\om) -
 \min_{I''_n} B(\cdot,\om) \Big), \frac{x-u}{v-u} \bigg) \, ,
\end{multline}
where the function $ \phi $ is defined by
\[
\phi(a,b,t) = \const(a,b) \int_0^{\min(a,b)} (a-y) (b-y) \exp \bigg( -\frac{
(a-y)^2 }{ 2t } - \frac{ (b-y)^2 }{ 2(1-t) } \bigg) \, \D y
\]
for $ t \in (0,1) $ and $ a,b > 0 $; the normalizing constant, $ \const(a,b)
$, ensures that $ \int_0^1 \phi(a,b,t) \, \D t = 1 $. (For $ \om \in
\Om\setminus C_n $ the formula is similar.) The formula follows easily from
the description of the conditional distribution given in the proof of Lemma
\ref{4.3}, the Brownian scaling, and the well-known joint distribution of the
minimizer $ T $ and the minimum $ Y = B_T $ of a Brownian path on $ [0,1] $
conditioned by $ B_0 = a $, $ B_1 = b $. Namely, the conditional density of $
(T,Y) $ is
\begin{equation}\label{4.5}
(t,y) \mapsto \sqrt{ \frac 2 \pi } \frac{ (a-y)(b-y) }{ (t-t^2)^{3/2} } \exp
\bigg( \frac{ (a-b)^2 }{ 2 } - \frac{ (a-y)^2 }{ 2t } - \frac{ (b-y)^2 }{
2(1-t) } \bigg)
\end{equation}
for $ 0 < t < 1 $, $ -\infty < y < \min(a,b) $.

We need the (unconditional) distribution of the random variable $ g_n(x,\cdot)
$ in order to check \ref{3.6}(c); the distribution should not concentrate near
the origin. However, the infimum of $ \phi(a,b,t) $ over all $ t \in (0,1) $
vanishes (unless $ a = b $). We restrict ourselves to a subinterval, say, the
inner half $ [1/4, 3/4] $ of $ [0,1] $; clearly, $ \phi(a,b,t) \ge \psi(a,b) >
0 $ for $ 1/4 \le t \le 3/4 $, therefore
\begin{equation}\label{4.6}
\Pr{ 0 < g_n(x,\cdot) < \eps } \le \xi(\eps) \, , \quad \xi(\eps) \to 0
\text{ as } \eps \to 0 
\end{equation}
for some $ \xi(\cdot) $ (not depending on $ n $ and $ x $), provided that $ x
$ belongs to the inner half of $ I'_n $ or $ I''_n $. (In fact we get much
more, namely, $ \Pr{ 0 < \operatorname{Length}(I_n) \cdot g_n(x,\cdot) < \eps
} \le \xi(\eps) $.)

\begin{theorem}\label{4.7}
The \valued{\DCS(0,1)} random variable $ \om \mapsto M_\om $ is distributed
like an unordered uniform infinite sample.
\end{theorem}

\begin{proof}
We will prove that the random variables $ X_n $ introduced in the proof of
Lemma \ref{4.1} satisfy the conditions of Lemma \ref{2.1}. First, we note that
almost every $ x \in (0,1) $ belongs to the inner half of $ I_n $ for
infinitely many $ n $. Let $ x $ be such a number; we will prove that $ \sum
f_n(x,\cdot) = \infty $ a.s.

We take $ n_1 < n_2 < \dots $ such that $ x $ belongs to the inner half of $
I'_{n_k} $ or $ I''_{n_k} $ for each $ k $, and define random variables $ Y_k,
Z_k $ by
\[
Y_k = f_{n_k}(x,\cdot) \, , \quad Z_k = g_{n_k}(x,\cdot) \, ,
\]
where $ f_n, g_n $ are the conditional densities introduced in the proof of
Lemma \ref{4.3}. (They are continuous in $ x $.) The relation $ f_n(x,\cdot) =
\cE{ g_n(x,\cdot) }{ \Ec_{n-1} } $, noted there, shows that $ Y_k = \cE{ Z_k
}{ X_1,\dots,X_{n_k-1} } $ which gives us \ref{3.6}(b). Condition \ref{3.6}(c)
follows from \eqref{4.6}. Taking into account Remark \ref{3.8} it remains to
prove that $ \{Y_k\ne0\} = \{Z_k\ne0\} = T^{-m_k} (A_k) $ for some $ m_1 < m_2
\dots $, some precompact sequence $ (A_k)_k $ such that $ \limsup_k \Pr{A_k}
> 0 $, and some mixing $ T : \Om \to \Om $.

We define $ T $ on the probability space of \emph{two-sided} Brownian paths as
the Brownian scaling centered at $ x $,
\[
B_{x+2s} (T\om) = \sqrt2 \( B_{x+s}(\om) - B_{x/2}(\om) \) \quad \text{for $ s
\in \R $} \, ;
\]
it is well-known to be mixing. Recalling the events $ C_n $ introduced in the
proof of Lemma \ref{4.3} we see that $ \{Y_k\ne0\} = \{Z_k\ne0\} = C_{n_k} $
for all $ k $ such that $ x \in I'_{n_k} $. (Other $ k $ satisfy $ x \in
I''_{n_k} $ and $ \{Y_k\ne0\} = \{Z_k\ne0\} = \Om \setminus C_{n_k} $; they are
left to the reader.) We have
\[
C_{n_k} = \Big\{ \min_{I'_{n_k}} B > \min_{I''_{n_k}} B \Big\} \, ,
\]
thus, $ C_{n_k} = T^{-m_k} (A_k) $ where $ m_k $ are such that the length of $
I_{n_k} $ is $ 2^{-m_k} $, and $ A_k $ are defined by
\[
A_k = \Big\{ \min_{[a_k-1/2,a_k]} B > \min_{[a_k,a_k+1/2]} B \Big\} \, ,
\]
$ a_k \in (x,x+1/2) $ being such that
\[
I_{n_k} = [ x - 2^{-m_k} (x-a_k+\tfrac12), x + 2^{-m_k} (a_k+\tfrac12-x) ] \,
.
\]
Clearly, $ \Pr{A_k} = 1/2 $ for all $ k $. Precompactness of the sequence $
(A_k)_k $ is ensured by continuity of the map $ a \mapsto \{ \min_{[a-1/2,a]}
B > \min_{[a,a+1/2]} B \} $ from $ \R $ to the space of events.
\end{proof}

\section[]{\raggedright The alternative way}
\label{sect5}In this section I abandon (temporarily!) my principle (formulated before
Def.~\ref{1.3}) and try nonstandard probability spaces. Given a standard Borel
space $ B $ and an equivalence relation $ E \subset B \times B $, the quotient
set $ B/E $ is equipped with the \sif\ $ \F_{B/E} $ of all sets $ A \subset
B/E $ whose inverse images in $ B $ (w.r.t.\ the canonical projection $ B \to
B/E $) are measurable. Thus, $ B/E $ is a Borel space (nonstandard, in
general).

In order to avoid ambiguity, concepts of Sect.~\ref{sect1} will be called
`strong', while concepts of this section --- `weak'. For example, a map $ \Om
\to B/E $ is strongly measurable, if it is measurable according to \ref{1.3},
and weakly measurable, if it is a measurable map from $ (\Om,\F,P) $ to $
(B/E, \F_{B/E}) $ according to the usual definition. (Still, $ (\Om,\F,P) $ is
a standard probability space.) Another example: weak distributions on $ B/E $
are just probability measures on $ (B/E, \F_{B/E}) $. Strong distributions are
much less customary objects (recall \ref{1.5}).

A strongly measurable map $ \Om \to B/E $ evidently is weakly measurable. The
converse is wrong in general (since $ E $ need not be measurable). Maybe it
holds under some reasonable condition on $ E $; I do not know.

If strongly measurable $ f : \Om_1 \to B/E $, $ g : \Om_2 \to B/E $ are
strongly identically distributed, then evidently they are weakly identically
distributed. We get a map from strong distributions on $ B/E $ to weak
distributions on $ B/E $. Is it injective? Is it surjective? I do not know.

Theorem \ref{4.7} considers two strong \valued{\DCS(0,1)} random variables and
states that they are strongly identically distributed. Therefore they are
weakly identically distributed, which allows us to transfer the Hewitt-Savage
zero-one law from the infinite sample to the Brownian minimizers, as follows.

\begin{corollary}\label{5.1}
Let $ U_1, U_2, \dots $ be i.i.d.\ random variables uniform on $ (0,1) $;
random variables $ X_1, X_2, \dots $ be all the Brownian local minimizers on $
(0,1) $ (enumerated as in Sect.~\ref{sect4} or otherwise); and $ A \subset
(0,1)^\infty $ a Borel set invariant under permutations. Then
\[
\Pr{ (X_1, X_2, \dots) \in A } = \Pr{ (U_1, U_2, \dots) \in A } \in \{0,1\}
\, .
\]
\end{corollary}

\begin{question}\label{5.2}
Let two strong \valued{\DCS(0,1)} random variables be weakly identically
distributed. Does it follow that they are strongly identically distributed?
(See also \ref{5.11}.)
\end{question}

\begin{proposition}\label{5.3}
If two strong \valued{\R/\Q} random variables are weakly identically
distributed then they are strongly identically distributed.
\end{proposition}

The proof is given after Proposition \ref{5.10}. Of course, by $ \R / \Q $ I
mean reals modulo rationals, that is, $ \R / E $ where $ E = \{ (x,y) \in \R^2
: x-y \in \Q \} $.

\begin{corollary}
Let probability measures $ \mu, \nu $ on $ \R $ be absolutely continuous
(w.r.t.\ the Lebesgue measure). Then there exists a probability measure $ J $
on $ \R^2 $, whose marginals are $ \mu, \nu $, such that
\[
x-y \in \Q \quad \text{for \almost{J} all $ (x,y) $.}
\]
\end{corollary}

You may try to construct such $ J $ explicitly, say, when $ \mu $ is uniform
and $ \nu $ is exponential.

Given a standard Borel space $ B $, we introduce the algebra $ \A $ of subsets
of $ B \times B $ generated by all product sets $ U \times V $ where $ U,V
\subset B $ are Borel sets. That is, elements of $ \A $ are of the form $ U_1
\!\times\! V_1 \,\cup\, \dots \,\cup\, U_n \!\times\! V_n $. The following
lemma is a slight modification of the well-known `marriage lemma'. By a
positive measure I mean a \valued{[0,\infty)} Borel measure (the measure of
the whole space is finite, and may vanish).

\begin{lemma}\label{5.5}
Let $ \mu,\nu $ be positive measures on $ B $, and $ W \in \A $. Then
\[
\sup_{m:m_1 \le\mu, m_2 \le\nu} m(W) = \inf_{U,V: W\subset\my U V}
\( \mu(U) + \nu(V) \) \, ;
\]
here $ m $ runs over positive measures on $ W $; $ U,V $ run over Borel
subsets of $ B $; and $ m_1, m_2 $ stand for the marginals of $ m $ (that is,
$ m_1(U) = m(U\times B) $ and $ m_2(V) = m(B\times V) $).
\end{lemma}

\begin{proof}
Clearly, $ \sup m(W) \le \inf \( \mu(U) + \nu(V) \) $ (since $ m(W) \le
m(U\times B) + m(B\times V) $); we have to prove that $ \sup m(W) \ge \inf \(
\mu(U) + \nu(V) \) $. First, we reduce the general case to the elementary case
of a finite set $ B $. To this end we take a finite partition $ B = B_1 \uplus
\dots \uplus B_n $ such that $ W $ is the union of $ B_k \times B_l $ (over
some pairs $ (k,l) $) and consider linear combinations of product measures $
(\mu \cdot \One_{B_k}) \times (\nu \cdot \One_{B_l}) $.

For a finite $ B $ we apply the usual duality argument in the
finite-dimensional space $ \R^B $:
\[
\sup m(W) = \inf_{f,g} \bigg( \int f \, \D\mu + \int g \, \D\nu \bigg)
\]
where the infimum is taken over all pairs of functions $ f,g : B \to
[0,\infty) $ such that $ f(x) + g(y) \ge 1 $ for all $ (x,y) \in W $. It
remains to prove that
\[
\inf \bigg( \int f \, \D\mu + \int g \, \D\nu \bigg) \ge \inf \( \mu(U)
+ \nu(V) \) \, .
\]
Introducing
\[
U_\theta = \{ x \in B : f(x) \ge \theta \} \, , \quad
V_\theta = \{ y \in B : g(y) \ge 1-\theta \}
\]
for $ \theta \in (0,1) $, we get $ W\subset\my{U_\theta}{V_\theta}
$ for each $ \theta $ and
\[
\int_0^1 \mu(U_\theta) \, \D\theta = \int f \, \D\mu \, , \quad
\int_0^1 \nu(V_\theta) \, \D\theta = \int g \, \D\nu \, ,
\]
therefore $ \int f \, \D\mu + \int g \, \D\nu \ge \inf_\theta \( \mu(U_\theta)
+ \nu(V_\theta) \) $.
\end{proof}

Here is a slight modification of a well-known result of Strassen
\cite[Sect.~6]{St} about measures with given marginals, concentrated on a
given \emph{closed} subset of a product space. A set of class $ \A_\de $ is,
by definition, a set of the form $ W_1 \cap W_2 \cap \dots $ where $ W_1, W_2,
\dots \in \A $ (and $ \A $ is introduced before \ref{5.5}).

\begin{lemma}\label{5.6}
Let $ \mu,\nu $ be positive measures on $ B $, and $ W \in \A_\de $. Then
\[
\sup_{m:m_1 \le\mu, m_2 \le\nu} m(W) = \inf_{U,V: W\subset\my U V} \( \mu(U) +
\nu(V) \) \, ;
\]
here $ m $ runs over positive measures on $ W $ and $ U,V $ run over Borel
subsets of $ B $ (and $ m_1, m_2 $ stand for the marginals of $ m $, as
before).
\end{lemma}

\begin{proof}
Once again, `$ \le $' is evident; we have to prove `$ \ge $'. We take $ W_1,
W_2, \dots \in \A $ such that $ W_n \downarrow W $ (that is, $ W_1 \supset
W_2 \supset \dots $ and $ W = W_1 \cap W_2 \cap \dots $). Lemma \ref{5.5}
applied to each $ W_n $ separately gives us measures $ m_n $ on $ W_n $
satisfying the restriction on marginals ($ (m_n)_1 \le \mu $, $ (m_n)_2 \le
\nu $) and such that
\[
m_n (W_n) + \frac1n \ge \inf_{U,V: W_n \subset\my U V} \( \mu(U) + \nu(V) \)
\ge \inf_{U,V: W\subset\my U V} \( \mu(U) + \nu(V) \) \, .
\]
The space of joinings, equipped with an appropriate topology, is a compact
metrizable space, and functions $ J \mapsto J(W) $ are continuous as long as $
W \in \A $; see the digression `The compact space of joinings' in
\cite[Sect.~4b]{Ts04}. This fact (and its proof) holds also for the space of
all positive measures $ m $ on $ B \times B $ satisfying $ m_1 \le \mu $, $
m_2 \le \nu $ (rather than $ m_1 = \mu $, $ m_2 = \nu $). Taking a
convergent subsequence $ m_{n_k} \to m $ we get
\[
m (W_n) \ge \inf_{U,V: W\subset\my U V} \( \mu(U) + \nu(V)
\)
\]
for all $ n $; however, $ m(W_n) \downarrow m(W) $.
\end{proof}

\begin{lemma}\label{5.7}
The following two conditions are equivalent for every $ W \subset B \times B
$:

(a) $ \inf \{ \mu(U) + \nu(V) : \my U V \supset W \} = 0 $;

(b) $ W\subset \my U V $ for some $ U,V $ such that $ \mu(U)=0 $, $ \nu(V)=0
$.
\end{lemma}

\begin{proof}
(b) \imp (a): trivial.

(a) \imp (b):
We take $ U_n, V_n $ such that $ W\subset \my{U_n}{V_n} $ for
each $ n $, and $ \sum ( \mu(U_n) + \nu(V_n) ) < \infty $. Then $ \mu (
\limsup U_n ) = 0 $ and $ \nu ( \limsup V_n ) = 0 $; here $ \limsup U_n $ is
the set of all $ x \in B $ such that $ x \in U_n $ for infinitely many $ n
$. It remains to note that $ W \subset \my{(\limsup U_n)}{(\limsup V_n)} $.
\end{proof}

We turn to measures with given marginals, concentrated on a given equivalence
relation $ E \subset B \times B $. By $ \F_E $ we denote the \sif\ of all
Borel sets $ A \subset B $ that are \saturated{E}, that is, $ (x,y) \in E \;
\& \; x 
\in A \imply y \in A $. If a measure $ m $ on $ B \times B $ is concentrated
on $ E $ (that is, $ ( B \times B ) \setminus E \subset A $, $ m(A) = 0 $ for
some Borel set $ A \subset B \times B $) then the marginal measures $ m_1, m_2
$ are equal on $ \F_E $ (that is, $ m_1(A) = m_2(A) $ for all $ A \in \F_E $),
since the symmetric difference between $ A \times B $ and $ B \times A $ is
contained in $ ( B \times B ) \setminus E $. By a \emph{nonzero} positive
measure I mean that the measure of the whole space does not vanish.

\begin{lemma}\label{5.8}
The following two conditions on $ E $ are equivalent:

(a) for every pair $ (\mu,\nu) $ of probability measures on $ B $ equal on $
\F_E $ there exists a probability measure $ m $ concentrated on $ E $ such
that $ m_1 = \mu $, $ m_2 = \nu $.

(b) for every pair $ (\mu,\nu) $ of nonzero positive measures on $ B $ equal
on $ \F_E $ there exists a nonzero positive measure $ m $ concentrated on $ E
$ such that $ m_1 \le \mu $, $ m_2 \le \nu $.
\end{lemma}

\begin{proof}
(a) \imp (b):
We note that $ \mu(B) = \nu(B) $, apply (a) to $ (1/\mu(B)) \mu $ and $
(1/\nu(B)) \nu $ and use $ \mu(B) m $.

(b) \imp (a):
We consider the set $ M $ of all positive measures $ m $ concentrated on $ E $
such that $ m_1 \le \mu $, $ m_2 \le \nu $. The set $ M $ contains a maximal
element $ m $, since $ M $ contains the limit of every increasing sequence of
elements of $ M $. We have to prove that $ m(B\times B) = 1 $. Assume the
contrary:
$ m(B\times B) < 1 $. The nonzero positive measures $ \mu - m_1 $, $ \nu - m_2
$ are equal on $ \F_E $. Item (b) gives us a nonzero positive measure $ \De m
$ concentrated on $ E $ such that $ (\De m)_1 \le \mu - m_1 $, $ (\De m)_2 \le
\nu - m_2 $. Thus, $ m + \De m $ belongs to $ M $, in contradiction to the
maximality of $ m $.
\end{proof}

\begin{remark}\label{5.9}
Let $ \mu,\nu $ be probability measures on $ B $ equal on $ \F_E $. Then the
following condition is sufficient for the existence of a probability measure $
m $ concentrated on $ E $ such that $ m_1 = \mu $, $ m_2 = \nu $:

(a) for every nonzero positive measures $ \mu_0,\nu_0 $ equal on $ \F_E $ and
satisfying $ \mu_0 \le \mu $, $ \nu_0 \le \nu $ there exists a nonzero
positive measure $ m' $ concentrated on $ E $ such that $ m'_1 \le \mu_0 $, $
m'_2 \le \nu_0 $.
\end{remark}

The proof is basically the same as the proof of `(b) \imp (a)' in Lemma
\ref{5.8}.

The saturation of a set $ A $ (w.r.t.\ a given equivalence relation $ E $) is,
by definition, $ \{ y \in B : \exists x \in A \; (x,y) \in E \} $. (It need
not be a Borel set even if $ A $ and $ E $ are Borel sets.) A set of class $
\A_{\de\si} $ is, by definition, a set of the form $ W_1 \cup W_2 \cup \dots $
where $ W_1, W_2, \dots \in \A_\de $ (and $ \A_\de $ is introduced before
\ref{5.6}).

\begin{proposition}\label{5.10}
Let $ B $ be a standard Borel space, $ E \subset B \times B $ an equivalence
relation of class $ \A_{\de\si} $ such that for every Borel set its saturation
is also a Borel set, and $ \mu,\nu $ probability measures on $ B $ equal on $
\F_E $. Then there exists a probability measure $ m $ concentrated on $ E $
such that $ m_1 = \mu $, $ m_2 = \nu $.
\end{proposition}

\begin{proof}
Assume that the sufficient Condition \ref{5.9}(a) is violated for some nonzero
$ \mu_0 \le \mu $, $ \nu_0 \le \nu $ equal on $ \F_E $. We take $ W_n \in
\A_\de $ such that $ E = W_1 \cup W_2 \cup \dots $ and note that each $ W_n $
violates the condition, that is, $ m=0 $ is the only positive $ m $
concentrated on $ W_n $ such that $ m_1 \le \mu_0 $, $ m_2 \le \nu_0 $. We
apply Lemma \ref{5.6} to $ \mu_0, \nu_0, W_n $; the supremum vanishes,
therefore the infimum vanishes. Lemma \ref{5.7} gives us $ U_n, V_n $ such
that $ \mu_0 (U_n) = 0 $, $ \nu_0 (V_n) = 0 $ and $ W_n \subset \my{U_n}{V_n}
$. Taking $ U = U_1 \cup U_2 \cup \dots $ and $ V = V_1 \cup V_2 \cup \dots $
we get
\[
\mu_0 (U) = 0 \, , \quad \nu_0 (V) = 0 \, , \quad E \subset \my U V \, .
\]
The latter means that a point of $ B \setminus U $ is never equivalent to a
point of $ B \setminus V $, that is, the saturation $ A $ of the set $ B
\setminus U $ is a subset of $ V $. We have $ A \in \F_E $ (the saturation of
a Borel set is Borel, as assumed), therefore $ \mu_0(A) = \nu_0(A) \le
\nu_0(V) = 0 $ and $ \mu_0 (B\setminus U) \le \mu_0(A) = 0 $, in contradiction
to the fact that $ \mu_0(B\setminus U) = \mu_0(B) > 0 $.
\end{proof}

Proposition \ref{5.3} is basically a special case of Proposition
\ref{5.10}. The equivalence relation $ E = \{ (x,y) \in \R^2 : x-y \in \Q \} $
belongs to the class $ F_\si $ (which means, the union of a sequence of closed
sets), therefore, to the class $ \A_{\de\si} $ (since every closed set belongs
to $ \A_\de $). The saturation $ A + \Q $ of any Borel set $ A $ is Borel
(since $ A + q $ is, for each $ q \in Q $). We have two strong \valued{\R/\Q}
random variables that are weakly identically distributed. They arise from two
\valued{\R} random variables whose distributions $ \mu, \nu $ are equal on $
\F_E $. Proposition \ref{5.10} gives us $ m $ concentrated on $ E $ whose
marginals are $ \mu, \nu $. This $ m $ is a joining between $ (\R,\mu) $ and $
(\R,\nu) $. It remains to lift the joining to the probability spaces, the
domains of our random variables, which is easy to do by means of the
conditional measures on these spaces.

In contrast, the equivalence relation \eqref{1.2} is of the class $ F_{\si\de}
$, therefore, $ \A_{\de\si\de} $.

\begin{question}\label{5.11}
Find a generalization of Proposition \ref{5.10} to equivalence relations of
the class $ \A_{\de\si\de} $, applicable to \eqref{1.2}. Is it possible?
(See also \ref{5.2}.)
\end{question}

\bigskip
\filbreak
{
\small
\begin{sc}
\parindent=0pt\baselineskip=12pt
\parbox{4in}{
Boris Tsirelson\\
School of Mathematics\\
Tel Aviv University\\
Tel Aviv 69978, Israel
\smallskip
\par\quad\href{mailto:tsirel@post.tau.ac.il}{\tt
 mailto:tsirel@post.tau.ac.il}
\par\quad\href{http://www.tau.ac.il/~tsirel/}{\tt
 http://www.tau.ac.il/\textasciitilde tsirel/}
}

\end{sc}
}
\filbreak

\end{document}